\documentclass[graybox]{svmult}

\usepackage{type1cm}        
\usepackage{makeidx}         
\usepackage{graphicx}        
\usepackage{multicol}        
\usepackage[bottom]{footmisc}

\usepackage{newtxtext}       %
\usepackage{newtxmath}       

\usepackage[caption=false]{subfig}
\usepackage{amsmath}
\usepackage{tikz}
\usetikzlibrary{patterns}
\usepackage{pgfplots}
\usepackage{pgfplotstable}
\usetikzlibrary{positioning}
\usetikzlibrary{calc}
\usetikzlibrary{decorations.pathreplacing}
\usepgflibrary{arrows}

\usepackage{soul}
\sethlcolor{pink}

\makeindex             

\begin{document}
\title*{A novel approach to fluid-structure interaction simulations involving large translation and contact}
\titlerunning{The phantom domain deformation mesh update method}
\author{Daniel Hilger, Norbert Hosters, Fabian Key, Stefanie Elgeti and Marek Behr}
\authorrunning{D. Hilger, N. Hosters, F. Key, S. Elgeti and M. Behr}
\institute{D. Hilger, N. Hosters, F. Key, S. Elgeti and M. Behr \at Chair for computational Analysis of technical Systems RWTH Aachen University, \mbox{Schinkelstra\ss{}e 2, 52062 Aachen}, {\{hilger\}\{hosters\}\{key\}\{elgeti\}\{behr\}@cats.rwth-aachen.de}}
%
%
\maketitle
\abstract*{In this work, we present a novel method for the mesh update in flow
problems with moving boundaries, the phantom domain deformation mesh update
method (PD-DMUM). The PD-DMUM is designed to avoid remeshing; even in the
event of large, unidirectional displacements of boundaries. The method combines
the concept of two mesh adaptation approaches: (1) The virtual ring shear-slip mesh
updatemethod (VR-SSMUM); and (2) the elastic mesh update method (EMUM). As
in the VR-SSMUM, the PD-DMUMextends the fluid domain by a phantom domain;
the PD-DMUM can thus locally adapt the element density. Combined with the
EMUM, the PD-DMUMallows the consideration of arbitrary boundary movements.
In this work, we apply the PD-DMUM in two test cases. Within the first test case,
we validate the PD-DMUM in a 2D Poiseuille flow on a moving background mesh.
Subsequently the fluid-structure interaction (FSI) problem in the second test case
serves as a proof of concept. More, we stress the advantages of the novel method
with regard to conventional mesh update approaches.}

\abstract{In this work, we present a novel method for the mesh update in flow
problems with moving boundaries, the phantom domain deformation mesh update
method (PD-DMUM). The PD-DMUM is designed to avoid remeshing; even in the
event of large, unidirectional displacements of boundaries. The method combines
the concept of two mesh adaptation approaches: (1) The virtual ring shear-slip mesh
update method (VR-SSMUM); and (2) the elastic mesh update method (EMUM). As
in the VR-SSMUM, the PD-DMUM extends the fluid domain by a phantom domain;
the PD-DMUM can thus locally adapt the element density. Combined with the
EMUM, the PD-DMUM allows the consideration of arbitrary boundary movements.
In this work, we apply the PD-DMUM in two test cases. Within the first test case,
we validate the PD-DMUM in a 2D Poiseuille flow on a moving background mesh.
Subsequently the fluid-structure interaction (FSI) problem in the second test case
serves as a proof of concept. More, we stress the advantages of the novel method
with regard to conventional mesh update approaches.}

\section{Introduction}\label{sec:1}
Many flow phenomena in technical processes, e.g., flows in liquid storage tanks, valve and piston flows, and in general all fluid-structure interaction problems involve moving boundaries.
The moving boundaries cause topological changes of the fluid domain which are important to consider when solving the flow problem.\\
\\
The changes of the fluid domain can be described either implicitly or explicitly \cite{ElgetiSauerland2016}.
In the implicit description--also called interface capturing--the boundary deformations are recorded on a fixed background mesh.
This strategy has the advantage that complex topology changes, e.g., breaking waves, can be resolved easily.
Yet, the treatment of discontinuities, conservation of mass, and the imposition of boundary conditions are still challenging.
Examples of interface-capturing methods are the volume-of-fluid method \cite{HirtNichols1981} or the level-set method \cite{OsherSethian1988}.
In the explicit description--called interface tracking--the domain deformations are described directly through the movements of its boundaries.
The mesh is restricted to the fluid domain and conforms with its boundaries.
This ensures an accurate approximation of the fluid interface and allows the imposition of boundary conditions along the moving boundary.
However, every time the topology of the domain is changed, the mesh must be adapted accordingly.\\
\\
The straightforward approach to incorporate the domain deformation is remeshing, but since remeshing is always connected to a projection of the solution between the old and the new mesh configurations, it should be avoided if possible \cite{JohnsonTezduyar1994}.
As an alternative to remeshing, mesh update methods can be used, where the current mesh is adapted to the changes of the domain.\\
\\
Mesh update methods can be categorized into two groups: (1) Methods in which the position of the mesh nodes are updated according to a predefined deformation rule, and (2) those where the mesh update is described by an additional set of equations \cite{Wall1999}.
In order to implement mesh update methods based on a predefined deformation rule, the changes of the fluid must be known in advance.
If this is not the case, the mesh update must be described by an additional equation.
One of the most commonly used methods of this type is the elastic mesh update method (EMUM) introduced in \cite{JohnsonTezduyar1994}.
Therein, the mesh is treated as an elastic body that deforms according to the motion of its boundaries.
For boundary movements that result in strongly constricted or expanded parts of the initial mesh, as it happens for example in valve flows, a mesh update will not provide a satisfactory solution.
This is because the existing mesh cells are either heavily squished or stretched.
In this case, remeshing of the fluid domain becomes inevitable..\\
\\
In order to avoid the need for remeshing, we propose a new mesh deformation method for large unidirectional mesh movements on boundary conforming meshes.
Therefore, we combine the EMUM and the  recently introduced virtual ring shear-slip mesh update method (VR-SSMUM) \cite{KeyPauliElgeti2018}.
The basic idea is here to perform the mesh update by means of the EMUM, but allow additional mesh cells to enter or exit the fluid domain.
Thus, the squeezing and the stretching of mesh cells is prevented by the possibility to increase or decrease the local number of finite elements (FE).
The new method is employed in conjunction with the deforming-spatial-domain/stabilized space-time (DSD/SST) approach \cite{TezduyarBehrLiou1992}, which is used to solve the flow problem on the changing domain.\\
\\
The structure of this paper is as follows:
In Section 2, we provide the governing equations of the flow problems we want to consider in the scope of this work. Further, we briefly summarize the DSD/SST method and the EMUM.
The concept and the implementation of the new mesh update method are explained in Section 3.
In Section 4, the validation and testing of the mesh update method is discussed by means of two test cases.

\section{Governing equations of fluid dynamics}
\label{sec:2}
The proposed mesh update method is developed specifically for flow problems with boundary conforming meshes involving large unidirectional boundary movements.
In this section, we present the governing equations of the flow problems examined within this work and further, we give a brief summary on the numerical methods employed to solve them.
\subsection{Governing equations of fluid dynamics}
Consider an incompressible fluid covering the deformable fluid domain $\Omega_t^f \subset  \mathrm{R}^{n_{sd}}$, with $n_{sd}$ indicating the number of spatial dimensions. At every time instant $t\in[0,T]$, the fluid's unknown velocity  $\mathbf{u}(\mathbf{x},t)$ and pressure $p(\mathbf{x},t)$ are governed by the Navier-Stokes equations for incompressible fluids:
\begin{subequations}
\begin{align}
\rho^f \left(\frac{\partial \mathbf{u}^f}{\partial t} \,+\, \mathbf{u}^f \cdot \boldsymbol{\nabla}  \mathbf{u}^f \,-\, \mathbf{f}^f \right)\,-\, \boldsymbol{\nabla}\cdot \boldsymbol{\sigma}^f\,=\, \mathbf{0} &\qquad\text{on} ~\Omega_t^f, \forall t \in \left(0,T\right),\\
\boldsymbol{\nabla} \cdot  \mathbf{u}^f\,=\, 0 &\qquad\text{on}~ \Omega_t^f, \forall t \in \left(0,T\right),
\end{align}
\label{Eq:NS}
\end{subequations}
with $\rho^f$ denoting the fluid density and $\mathbf{f}^f$ representing all external body forces per unit mass.
For Newtonian fluids, the stress tensor $\boldsymbol{\sigma}^f$ is defined as
\begin{equation}
\boldsymbol{\sigma}^f\,=\,-p^f \mathbf{I}\,+\,2\rho^f\nu^f \boldsymbol{\varepsilon}^f(\mathbf{u}^f),
\end{equation}
with
\begin{equation}
\boldsymbol{\varepsilon}^f(\mathbf{u}^f)\, =\, \frac{1}{2} \left( \boldsymbol{\nabla} \mathbf{u}^f + \left(\boldsymbol{\nabla} \mathbf{u}^f \right)^T \right),
\end{equation}
where $\nu^f$ denotes the dynamic viscosity. A well-posed system is obtained when boundary conditions are imposed on the external boundary  $\Gamma^{f}_{t}$. Here, we distinguish between Dirichlet and Neumann boundary conditions given by:
\begin{subequations}
	\begin{align}
	\mathbf{u}^f \,=\,\mathbf{g}^f &\qquad\text{on}~\Gamma^{f}_{t,g},\\
	\mathbf{n}^f\cdot \boldsymbol{\sigma}^f \,=\,\mathbf{h}^f &\qquad\text{on}~\Gamma^{f}_{t,h},
	\end{align}
\end{subequations}
where $\mathbf{g}^f$ and $\mathbf{h}^f$ prescribe the velocity and stress values on complementary subsets of $\Gamma^{f}_{t}$.
With regard to deformation of the fluid domain $\Omega_t^f$ in time, the DSD/SST method is applied to solve the Navier-Stokes equations.
\subsection{Deforming-spatial-domain/stabilized space-time method}\label{sec:DSDSST}
\begin{figure}[h]
	\centering
	\includegraphics{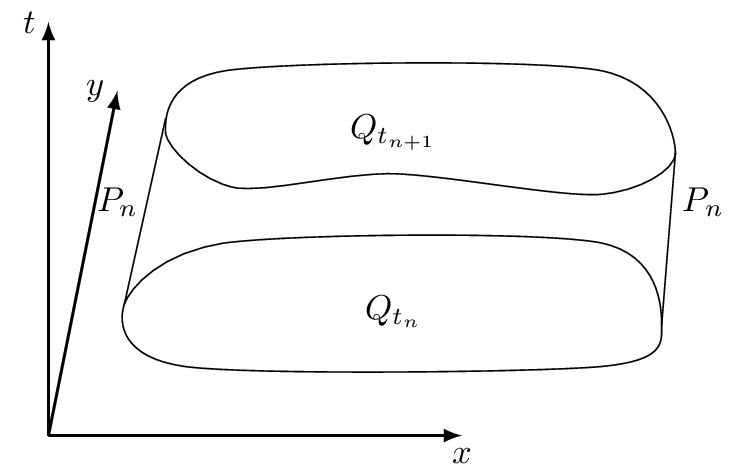}
	\caption{Space-time slab.}
	\label{fig:SpaceTimeSlab}
\end{figure}

The DSD/SST method is a space-time-based finite-element (FE) method, i.e., a FE discretization is applied to space and time.
It was first applied to flow problems with moving boundaries in \cite{TezduyarBehrLiou1992,TezduyarEtAl1992}.\\
\\
The advantage of the DSD/SST method is, that the variational form of the governing equations implicitly incorporates the deformations of the domain.
In order to construct the interpolation and weighting function spaces used in the variational formulation of the problem, the time interval $(0,T)$ is split into $N$ subintervals $I_n=\left[t_n,t_{n+1}\right]$, where $t_n$  and $t_{n+1}$ belong to an ordered series of time levels.
Thus, the space-time continuum is divided into $N$ space-time slabs $Q_n$ as depicted in Figure \ref{fig:SpaceTimeSlab}, bounded by the spatial configurations $\Omega_t$ at time $t_n$ and $t_{n+1}$, and $P_n$ describing the course of the spatial boundary $\Gamma^{f}_t$ as $t$ traverses $I_n$.
The boundary $P_n$ can be decomposed into two complementary subsets $(P_n)_g$ and $(P_n)_h$, representing the Dirichlet and Neumann boundary conditions of $\Gamma^f_{t}~\forall t \in I_n$.
The space-time slabs are weakly coupled along their interfaces using jump terms.
For the spatial approximation $\Omega^{f}_{t,h}$ of the domain $\Omega^f_t$, the following finite element trial and weighting function spaces are constructed:
\begin{subequations}
	\begin{align}
	\mathcal{H}^{1h}(Q_n)&:= \left\{ \mathbf{w}^h \in \mathcal{H}^1 \left(Q_n\right) \left| \mathbf{w}^h_T\right| \text{is a first-order polynominal }\forall T\in \mathcal{T}^h \right\},\\
	\mathcal{S}^h_u &:= \left\{ \mathbf{u}^h | \mathbf{u}^h \in \left[\mathcal{H}^{1h}\left(Q_n\right) \right]^{nsd}, \mathbf{u}^h =\mathbf{g} ~\text{on} ~\left(P_n\right)_g \right\},\\
	\mathcal{V}^h &:= \left\{ \mathbf{w}^h | \mathbf{w}^h \in \left[\mathcal{H}^{1h} \left(Q_n\right) \right]^{nsd}, \mathbf{w}^h= \mathbf{0} ~\text{on} ~\left(P_n\right)_g \right\},\\
	\mathcal{S}^h_p &=  \mathcal{V}^h_p :=\left\{ {q}^h | {q}^h \in \mathcal{H}^{1h}\left(Q_n\right)  \right\}.
	\end{align}
\end{subequations}
The interpolation functions are globally continuous in space, but discontinuous in time.
Using the following notational convention,
\begin{subequations}
\begin{align}
\left(\mathbf{u}^h\right)^{\pm}_n\,=\,\underset{\epsilon \rightarrow 0}{lim}~  \mathbf{u}\left(t_n \pm \epsilon \right)\\
\int_{Q_n} \cdots \text{d} Q = \int_{I_n} \int_{\Omega_t} \cdots \text{d}\Omega \text{d}t,\\
\int_{(P_n)} \cdots \text{d} P = \int_{I_n} \int_{\Gamma_t} \cdots \text{d}\Gamma \text{d}t,
\end{align}
\end{subequations}
and following references \cite{TezduyarBehrLiou1992, HughesFrancaHulbert1989, PauliBehr2017}, the stabilized variational formulation of the Navier Stokes equations is obtained:
Given $\left(\mathbf{u}^h\right)_n^-$ with $\left(\mathbf{u}^h\right)_0^- = \mathbf{u}_0$, find $\mathbf{u}^h \in \mathcal{S}^h_\mathbf{u}$ and $p^h\in \mathcal{S}^h_p$ such that $\forall \mathbf{w}^h \in \mathcal{V}^h_\mathbf{u}$, $\forall q\in \mathcal{V}^h_P$:
\begin{align}
\int_{Q_n} \mathbf{w}^h \cdot \rho^f \left( \frac{\partial \mathbf{u}^h}{\partial t} + \mathbf{u} \cdot \boldsymbol{\nabla} \cdot \mathbf{u}^h - \mathbf{f}\right) \text{d}Q +
\int_{Q_n} \boldsymbol{\nabla}\mathbf{w}^h : \boldsymbol{\sigma} ( p^h, \mathbf{u}^h ) \text{d}Q \nonumber \\
+ \int_{Q_n} q^h \boldsymbol{\nabla}\cdot \mathbf{u}^h \text{d}Q  +
\int_{\Omega_n} \left(\mathbf{w}^h\right)^+_n \cdot \rho^f \left(\left(\mathbf{u}^h\right)^+_n - \left(\mathbf{u}^h\right)^-_n\right)\text{d}\Omega \nonumber \\
+ \sum_{e=1}^{n_{el}}\int_{Q^e_n} \frac{1}{\rho^f} \tau_{MOM}  \left[\rho^f \mathbf{u}^h \cdot \boldsymbol{\nabla} \mathbf{w}^h + \boldsymbol{\nabla} q^h\right] \nonumber \\
 \cdot \left[\rho^f \left(\frac{\partial \mathbf{u}^h}{\partial t} + \mathbf{u} \cdot \boldsymbol{\nabla} \cdot \mathbf{u}^h - \mathbf{f} \right) - \boldsymbol{\nabla}\cdot\boldsymbol{\sigma}(p^h, \mathbf{u}^h) \right]\text{d}\Omega \nonumber\\
+ \sum_{e=1}^{n_{el}}\int_{Q_n^e} \boldsymbol{\nabla} \cdot \mathbf{w}^h \rho^f \tau_{CONT} \boldsymbol{\nabla} \cdot \mathbf{u}^h \text{d}\Omega \nonumber\\
= \int_{\left(P_n\right)_h} \mathbf{w}^h \cdot \mathbf{h}^h\text{d}P. \label{eq:NS_WEAK}
\end{align}
In Equation \eqref{eq:NS_WEAK}, the first three terms and the last term directly result from the variational formulation of Equation \eqref{Eq:NS}, whereas the fourth term denotes the jump terms between the space-time slabs.
Terms five and six result from a Galerkin-Least Squares (GLS) stabilization applied to the Navier-Stokes equations.
The stabilization approach used within this work and the choice of the stabilization parameters $\tau_{CONT}$ and $\tau_{MOM}$ are described in detail in \cite{PauliBehr2017}.\\
\\
Though the DSD/SST method implicitly accounts for the domain deformations in one time slab, a deformation rule is needed to deform the FE mesh according to the boundary movements.
\subsection{Elastic mesh update method}\label{sec:EMUM}
One approach for the automatic mesh update in boundary conforming meshes is the elastic mesh update method (EMUM) introduced by \cite{JohnsonTezduyar1994}, where the mesh is understood as an elastic body occupying the bounded region \mbox{$\Omega^\# \subset \mathcal{R}^{n_{sd}}$} with boundary $\Gamma^\#$.
Thus, the deformation of the mesh is expressed in terms of the nodal displacements $\mathbf{d}^\#$ governed by the equilibrium equation of elasticity:
\begin{equation}
\boldsymbol{\nabla} \cdot \boldsymbol{\sigma}^\#\,=\,\mathbf{0},
\end{equation}
where $\boldsymbol{\sigma}^\#$ corresponds to the Cauchy stress tensor,
\begin{equation}
\boldsymbol{\sigma}^\# \,=\, \lambda \left( tr \boldsymbol{\epsilon}^\# \right) \mathbf{I}\,+\,2\mu\boldsymbol{\epsilon}^\# ~,\qquad
\boldsymbol{\epsilon}^\#\,=\, \frac{1}{2} \left( \boldsymbol{\nabla} \mathbf{d}^\# + \left(\boldsymbol{\nabla} \mathbf{d}^\# \right)^T \right).
\end{equation}
The imposition of Dirichlet and Neumann boundary conditions yields a well-posed problem for the mesh deformation:
\begin{align}
\mathbf{d}^\# \,=\,\mathbf{g}^\# &\qquad\text{on}~ \left(\Gamma\right)^\#_g,\\
\mathbf{n}\cdot \boldsymbol{\sigma}^\# \,=\,\mathbf{h}^\# &\qquad\text{on}~ \left(\Gamma\right)^\#_h,
\end{align}
where $\mathbf{g}^\#$ and $\mathbf{h}^\#$ prescribe the displacements and normal stresses on the mesh boundaries.\\
\\
The elasticity problem is solved with the Galerkin FE method and the resulting displacements are applied to the mesh nodes representing the upper mesh configuration of the current space-time slab.

\section{The phantom domain mesh deformation method}
\label{sec:3}
\begin{figure}[h]
	\centering
	\resizebox {\textwidth} {!}{
	\includegraphics{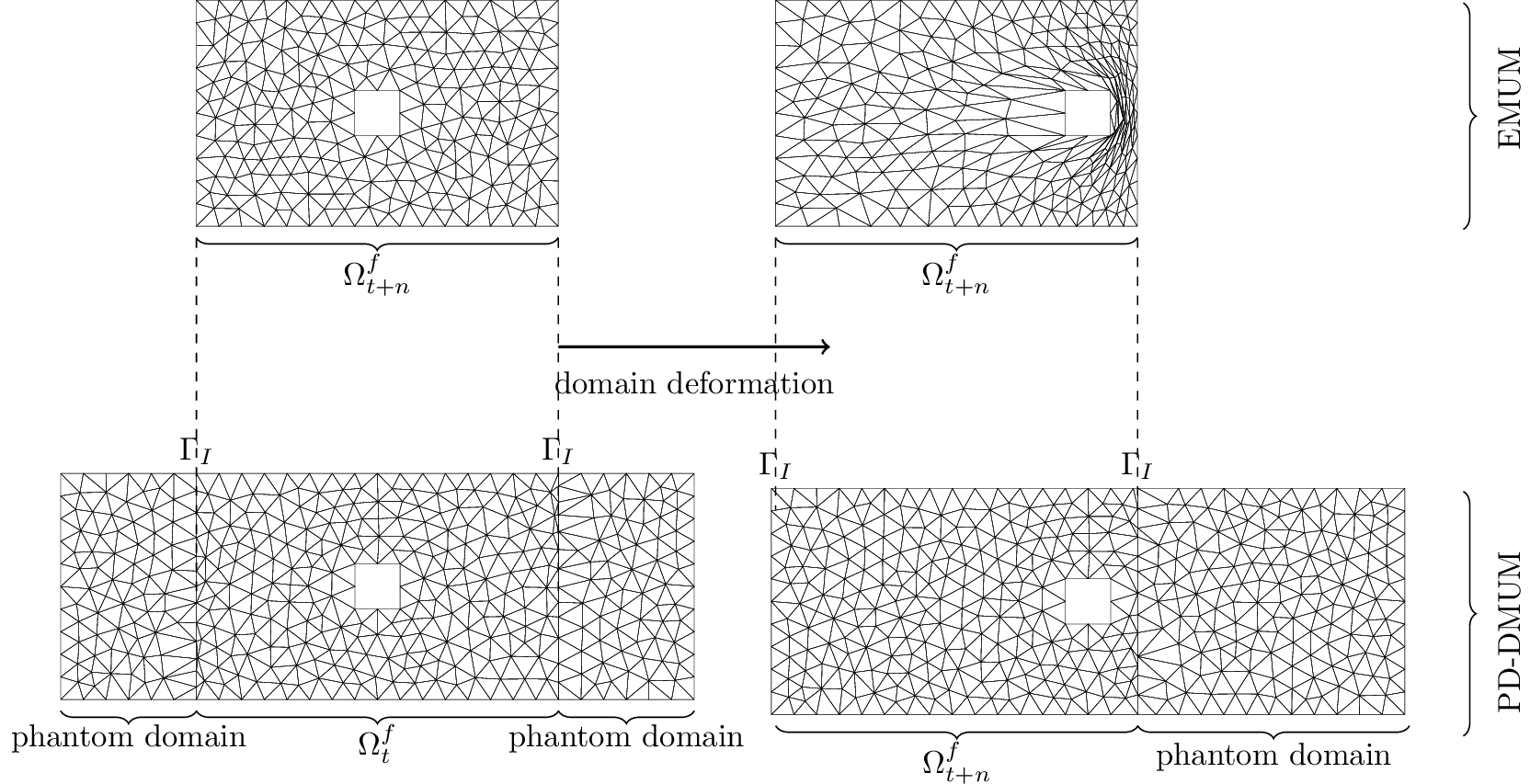}}
	\caption{Mesh deformation with PD-DMUM vs. EMUM}
	\label{fig:SquishedElements}
\end{figure}

\begin{figure}[h]
	\centering
	\resizebox {\textwidth} {!}{
		\includegraphics{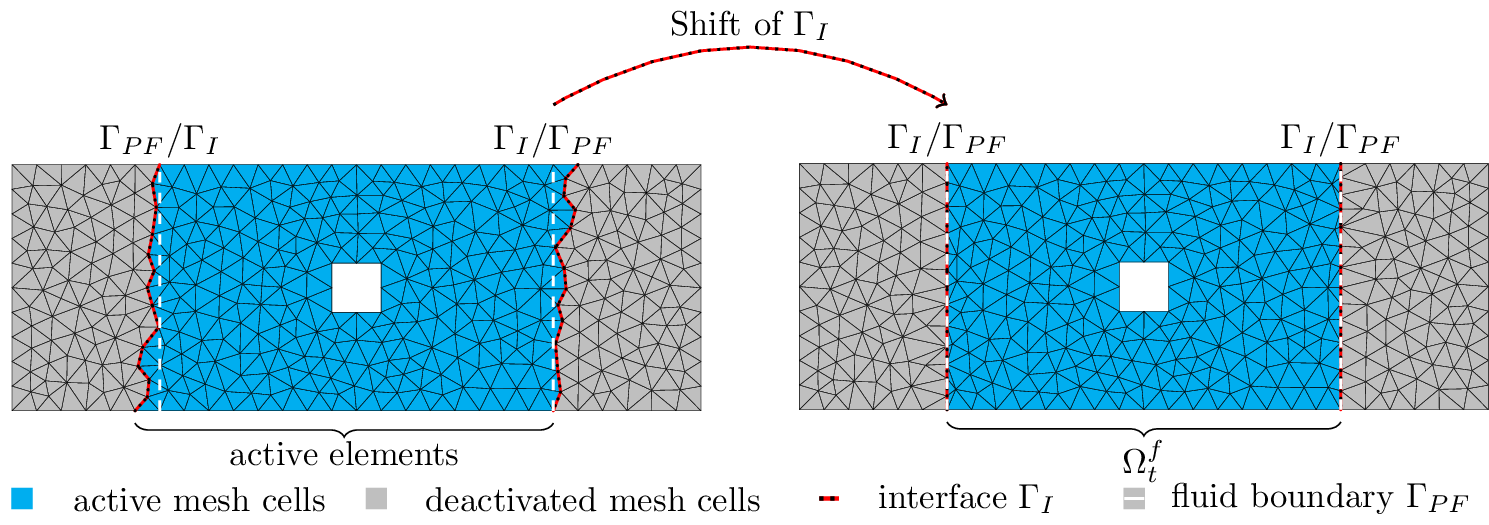}}
	\caption{Activity pattern on mesh with initial and corrected course of interface $\Gamma_I$.}
	\label{fig:ActivityBndry}
\end{figure}

\begin{figure}[ht!]
\centering
	\subfloat[\label{fig:Space-TimeUnshifted} uncorrected boundary $\Gamma_{I}$.]{
		\resizebox{0.45\textwidth}{!}{
			\includegraphics{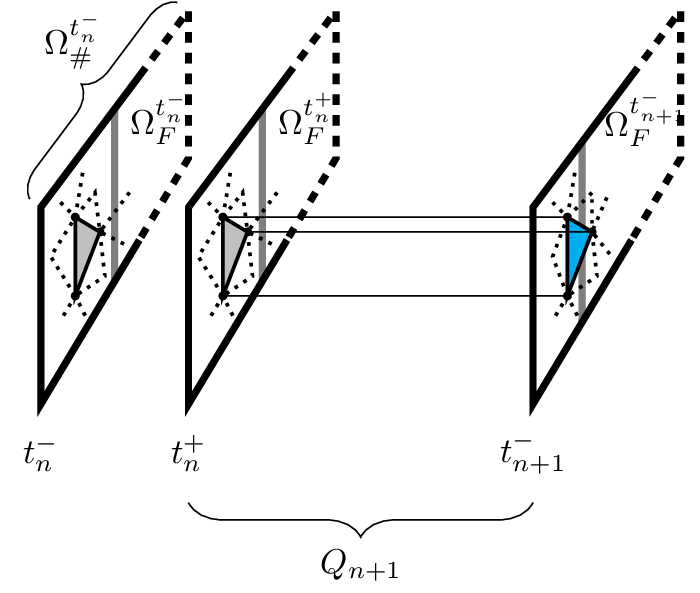}}
		}
	\subfloat[\label{fig:Space-TimeShifted} corrected boundary $\Gamma_{I}$.]{
		\resizebox{0.45\textwidth}{!}{
			\includegraphics{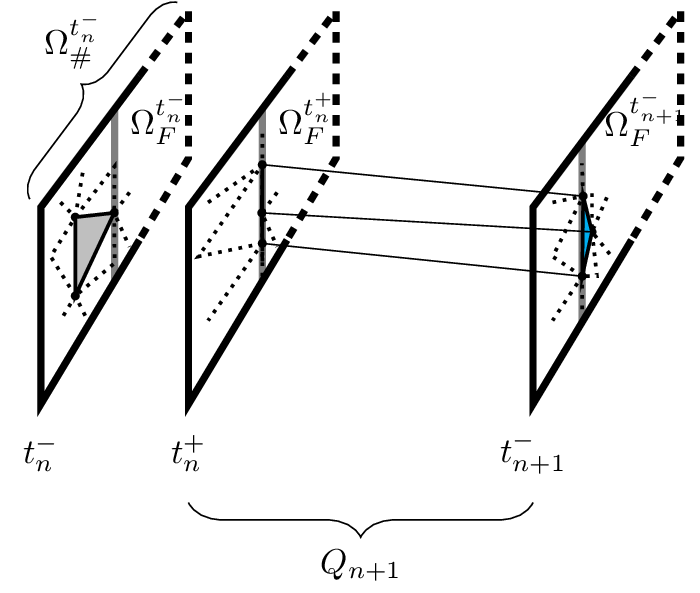}}
	}
	\caption{Shifting of $\Gamma_I$ for boundary conformity in space-time slab.}
	\label{fig:Space-TimeInterface}
\end{figure}

\begin{figure}[h]
	\centering
	\resizebox {\textwidth} {!}{
	\includegraphics{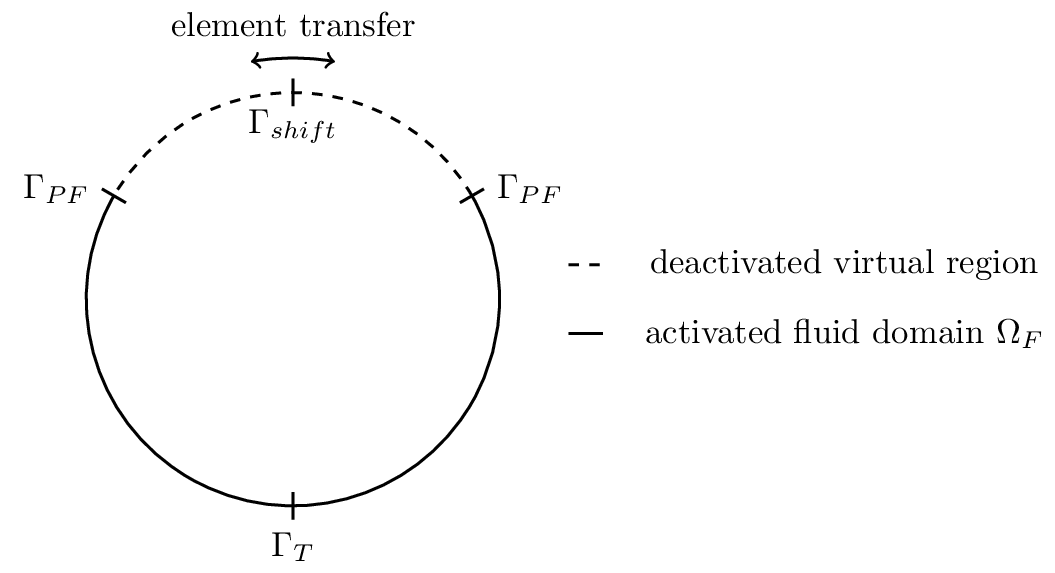}}
	\caption{Illustration of the virtual ring concept.}
	\label{fig:VirtualRing}
\end{figure}

The aim of the newly proposed method is to extend the usability of boundary-conforming meshes for deforming domains with large, unidirectional deformations.
The specific target are applications with large, unidirectional deformations (imagine an object sinking within a fluid or the flow through a valve).
So far, the fluid domain is enclosed within two types of boundaries: (1) deforming, and (2) fixed.
The deforming boundaries are handled in a standard interface tracking way, meaning that the boundary deforms according to its relevant deformation rule -- e.g., determined by the structure in an FSI context or a free-surface motion -- while the inner nodes adapt to this motion.
As depicted in \mbox{Figure \ref{fig:SquishedElements},} a predominantly unidirectional deformation, however, soon results in a situation where one side of the mesh contains very compressed elements, whereas the other side is comprised of very stretched elements.\\
\\
In our proposed method, this is counteracted via the implementation of a new boundary condition that allows mesh cells to exit and enter the fluid domain as needed.
The implementation of this boundary condition is based on the concept of the VR-SSMUM presented in \cite{KeyPauliElgeti2018}.
As with the VR-SSMUM, the mesh is extended by additional mesh cells.
As sketched in Figure \ref{fig:SquishedElements}, these cells are positioned in a phantom domain which is located outside of the fluid domain.
In the following we will therefore refer to this method as the phantom domain deformation mesh update method (PD-DMUM).\\
\\
Since not all mesh cells are positioned within the fluid domain, an activity pattern, as illustrated in Figure \ref{fig:ActivityBndry}, is used to determine which elements are used in the computation of the flow problem.
Here, elements that intersect with the fluid domain are considered as activated elements whereas the remaining elements are deactivated.
Activated and deactivated elements have a common interface $\Gamma_{I}$.
The interface is  a boundary of the fluid domain, which requires the definition of boundary values.
The boundary value prescribed at the element nodes of the interface is of a new boundary type.
The element nodes associated with the new boundary type have the special characteristic that they prescribe boundary values to the flow problem, but function as internal nodes in the mesh update method.
Consequently, the mesh of the phantom domain and the fluid domain are considered as one coherent mesh in the mesh deformation process.\\
\\
Now that the mesh is deformed according to the underlying deformation rule, elements from the phantom domain can slide across the prescribed fluid boundary $\Gamma_{PF}$ into the fluid domain or vice versa.
This changes the composition of elements that intersect with the fluid domain, so that the activity pattern of the elements must be re-determined.
In the space-time approach used here, one space-time slab is bounded by two different mesh configurations. This can lead to the situation shown in Figure \ref{fig:Space-TimeUnshifted}, where an element is located inside the fluid domain on the upper time level, yet positioned outside at the lower time level.
Therefore, we define here that the mesh configuration at the upper time level always determines which elements represent the fluid domain.
Based on the updated activity pattern, the new location of the interface $\Gamma_I$ is determined within the mesh.
The position of $\Gamma_I$ usually does not correspond to the position of the predefined fluid boundary $\Gamma_{PF}$.\\
\\
The boundary conformity of the mesh for $\Gamma_{PF}$ is now obtained by a closest point projection of all mesh nodes on $\Gamma_{I}$ to the prescribed contour of the fluid boundary.
It is important to note the special case of those elements which were not yet part of the fluid discretization in the previous time step, because these elements require a projection of the old solution onto the new boundary nodes.
This is necessary to calculate the jump terms in Equation \eqref{Eq:NS}.
This means that the new method does not require remeshing, yet the projection between two mesh configurations cannot be completely avoided.
However, the projection is limited to single elements when they enter the fluid domain.\\
\\
The sequence of the individual steps within the PD-DMUM can be summarized as follows:
\begin{enumerate}
\item Update mesh according to moving boundaries.
\item \textbf{Identify activated and deactivated elements.}
\item \textbf{Adapt the boundaries to the prescribed position of the fluid domain.}
\item \textbf{Set boundary values for the nodes on the redefined interface $\Gamma_I$.}
\item \textbf{Project the solution of the previous time step for all newly activated elements.}
\item Solve flow problem on active elements.
\end{enumerate}
In direct comparison with a conventional update strategy for boundary conforming meshes, such as the EMUM, steps (2)--(5) are those which are additionally required.\\
\\
Depending on the boundary movements, the PD-DMUM can be complemented with additional mesh update strategies.
In case of large unidirectional boundary movements, we can employ the concept of the virtual ring presented in \cite{KeyPauliElgeti2018}.
The objective of the virtual ring is to reduce the size of the phantom domain in the mesh update.
For this purpose, we connect the mesh along the outward facing boundaries of two oppositely positioned phantom domains.
This connection results in a coherent mesh, forming a virtual ring as illustrated in Figure \ref{fig:VirtualRing}.
The mesh update can now transfer elements between the connected phantom domains, while moving them along the virtual ring.
Consequentially, elements can exit the fluid domain on one side and re-enter the domain on the other side.
Therefore, the phantom domains can be reduced to a thin layer of elements.
The reduction of the phantom domains results in a significant decrease of computational cost for the mesh deformation problem.

\section{Computational results}
\label{sec:4}
The implementation of the PD-DMUM is applied to two test cases.
In a first step, we validate the mesh update method by examining its influence on the solution of a two dimensional Poiseuille flow.
In the second test case we show, by means of an example from the field of FSI, the advantages of the PD-DMUM.
\subsection{2D Poiseuille flow on moving background mesh}
\begin{figure}[h]
		\resizebox {\textwidth} {!}{
	\includegraphics{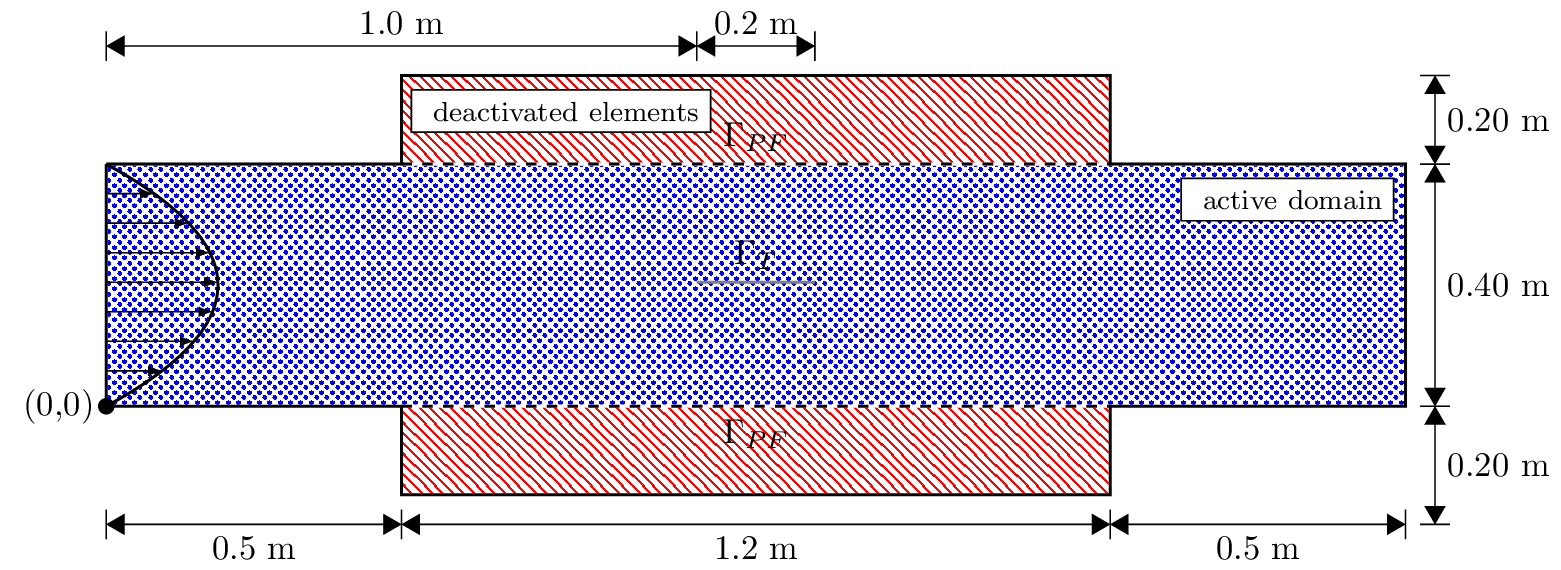}}
	\centering
	\caption{Tube geometry for Poiseuille flow.}
	\label{fig:PoiseuilleGeometry}
\end{figure}

\begin{figure}[h]
	\resizebox{\textwidth}{!}{
	\includegraphics{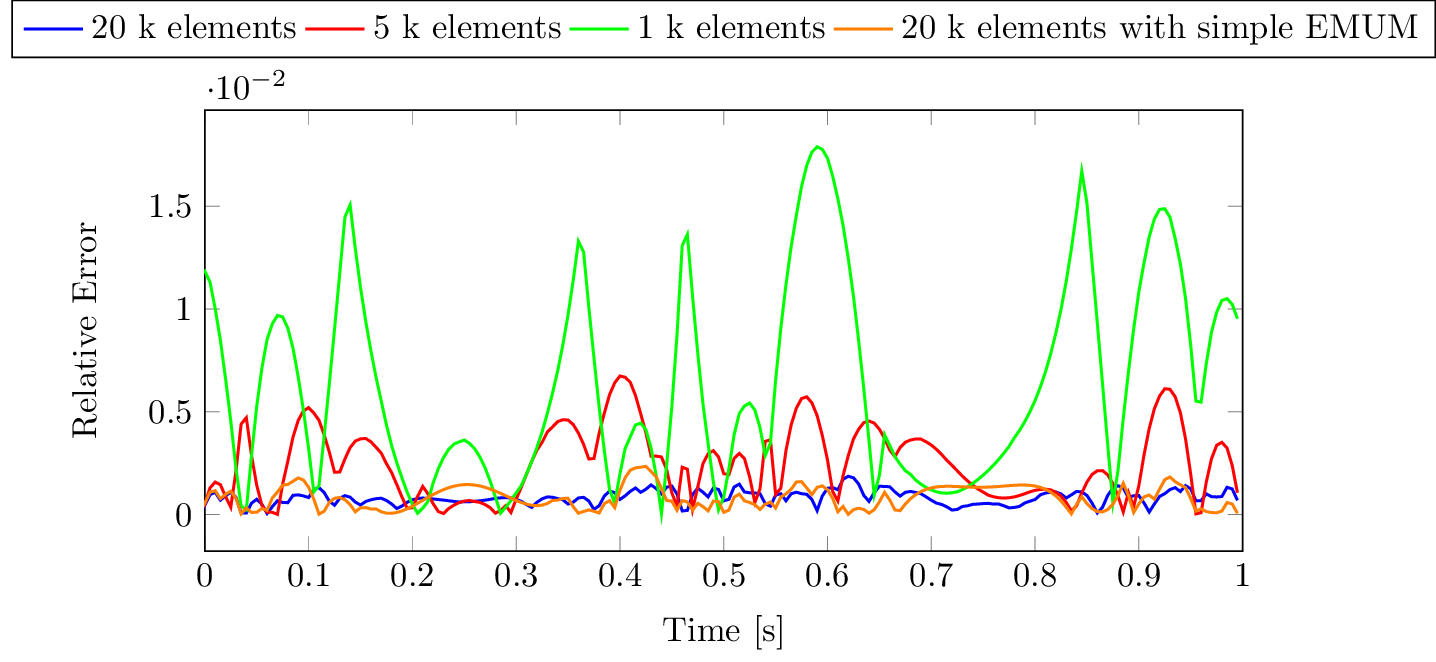}}
	\centering
	\caption{Relative error of velocity at probe position.}
	\label{fig:RelativErrorPoiseuille}
\end{figure}

\begin{figure}[h]
	\resizebox{\textwidth}{!}{
	\includegraphics{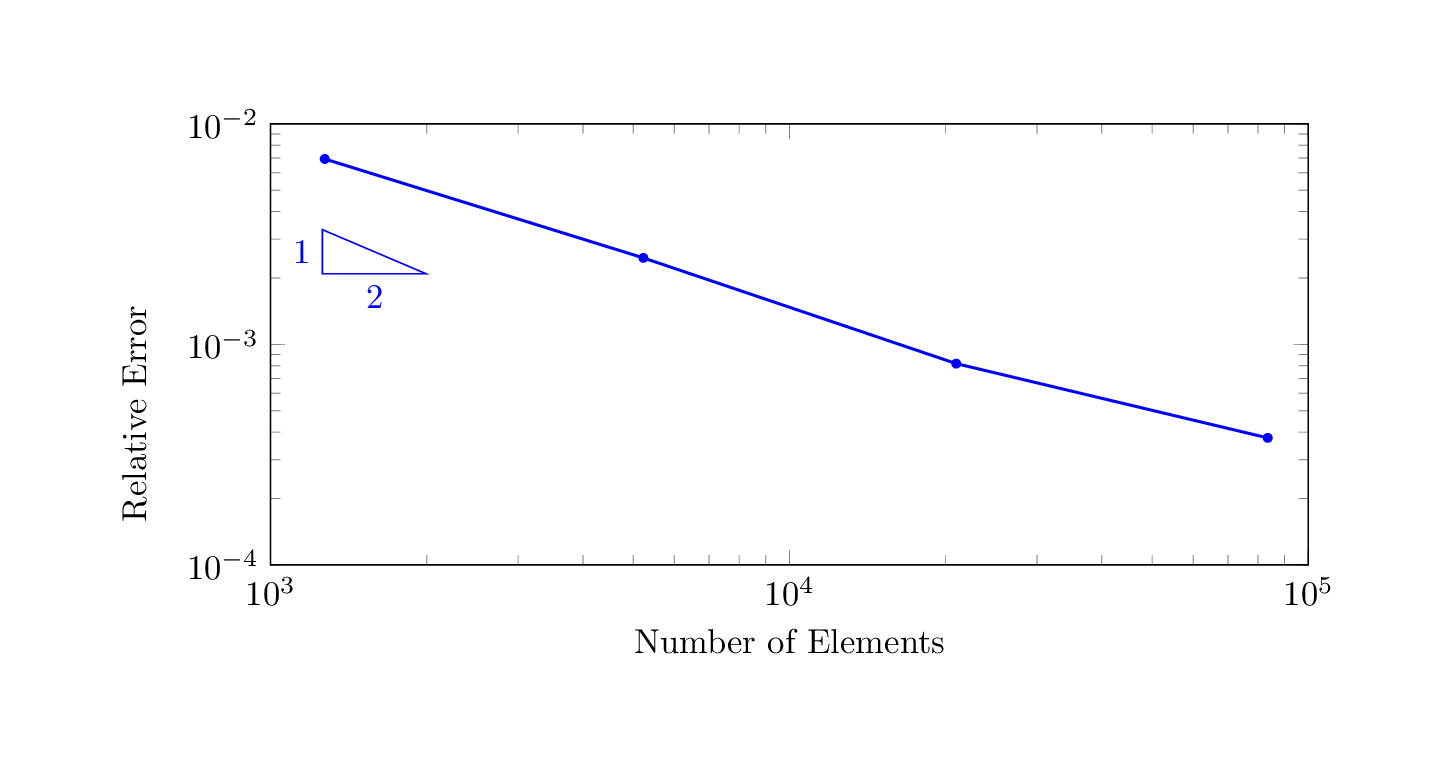}}
	\centering
	\caption{Average relative error for different mesh resolutions.}
	\label{fig:ConvergencePoiseuille}
\end{figure}

In the first test case we examine the influence of the PD-DMUM on a flow problem with a well-known solution.
For this purpose, we consider a two-dimensional Poiseuille flow in a tube.
The topology of the fluid domain remains unchanged, yet a predefined motion is applied to the underlying mesh.
The PD-DMUM is used to perform the mesh update, but should not affect the flow field within the tube.\\
\\
The geometric dimensions of the tube are chosen  according to Figure \ref{fig:PoiseuilleGeometry}.
In the middle of the domain, we position a mesh section $\Gamma_T$ by means of which the predefined mesh motion is imposed as a Dirichlet boundary condition.
The boundary $\Gamma_T$ has no physical impact with respect to the flow problem.
The additional phantom domains required within the PD-DMUM are positioned along the upper and lower boundary of the tube.
The material properties of the fluid are chosen according to Table \ref{tab:PoiseuilleFlow}.
\begin{table}[h]
	\begin{tabular}{ccc}
		\hline
		Parameter & Identifier  & Value  \\
		\hline
		density & $\rho$  &  $1.0$ [kg/m$^3$]\\
		viscosity & $\nu$  & $0.001$ [kg/m$\cdot $s] \\
		mean velocity& $U$  & $2.5$ [m/s]  \\
		\hline\\
	\end{tabular}
	\caption{Properties of fluid in 2D Poiseuille flow.}
	\label{tab:PoiseuilleFlow}
\end{table}
Regarding boundary conditions of the flow, we impose no-slip condition along the walls of the tube.
This also applies to the boundary section $\Gamma_{PF}$ at the interface between the phantom domain and the fluid domain.
A parabolic inflow profile for the velocity is given at the inlet of the tube:
\begin{equation}
\mathbf{u}(y) \,=\,\left( \frac{4Uy(H-y)}{H^2}, 0 \right).
\end{equation}
With respect to the mesh update, the position of the nodes at the inlet, the outlet, and the tube walls are fixed.
However, this does not apply to $\Gamma_{PF}$ and the remaining boundaries of the phantom domain, as these nodes should be able to move freely.
For the boundary $\Gamma_T$ we prescribe the following sinusoidal movement:
\begin{equation}
\mathbf{d}(t) \,=\,\left(0~,~ 0.1\cdot\text{sin}\left(\frac{2\,\pi\, t}{T} \right) \right).
\end{equation}
The mesh deformation is examined for a period of $T=8$[s]. The time step size is $\Delta t = 0.02$ [s].
Initially, a fully developed flow profile is already present in the pipe.\\
\\
The Poiseuille flow is computed on four mesh configurations with the PD-DMUM and for the purpose of comparison for one configuration by the EMUM.
For the comparison of the solutions we use the flow velocity.
The velocity is measured at a probe positioned at point $(1.1\,,\,0.2)$ inside the tube.
Together with the given analytical solution of the Poiseille flow, the relative error can be computed for the different mesh configurations.\\
\\
In a first step, the relative error of the computed velocity is evaluated for the probe position.
In Figure \ref{fig:RelativErrorPoiseuille} it can be observed that the relative error decreases as the mesh is refined.
The comparison between the solution of the EMUM and the PD-DMUM on similar grids shows that the relative error for the calculated velocity is of the same order of magnitude.
The fluctuations that can be observed for all computations can be explained by the linear interpolation of the parabolic velocity profile at the probe position.
In Figure \ref{fig:ConvergencePoiseuille}, we can observe that the numerical solution converges for the PD-DMUM towards the analytic solution of the Poiseuille problem.
Both, the convergence of the PD-DMUM and the comparable results to the EMUM for moderate mesh deformations indicate that the PD-DMUM provides a valid mesh update.
\subsection{Falling ring in a fluid-filled container}
\begin{figure}[!h]
		\resizebox {0.48\textwidth} {!}{
	\includegraphics{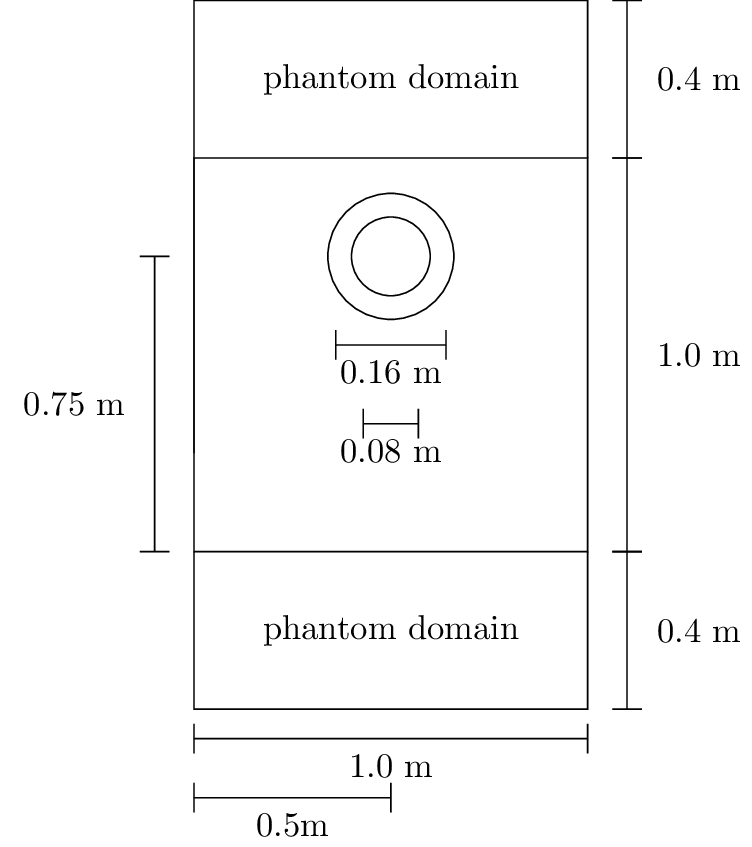}}
	\centering
	\caption{Geometry of container with ring.}
	\label{fig:SketchCylinder}
\end{figure}

The second test case is used to illustrate possible applications of the PD-DMUM.
For this purpose, we consider a fluid-structure interaction with large translational boundary movement.
More precisely, we simulate an elastic ring that falls inside a fluid-filled container until it hits the ground and rebounds.
Concerning the mesh deformation, this is a demanding process, since the number of mesh cells, which are initially positioned between the ring and the bottom, must be reduced to zero by the time of contact.
Using previous mesh update methods it is not possible to simulate this process on boundary conforming meshes without frequent remeshing of the fluid domain.\\
\\
The geometric dimensions of the container and the ring are chosen according to Figure \ref{fig:SketchCylinder}.
The ring is represented by a non-uniform rational B-spline (NURBS) \cite{PieglTiller1997} with 721 elements and second-order basis functions.
In total 13448 elements are used to discretize the fluid domain and the additional phantom domains.
In the flow problem no-slip conditions are prescribed along the walls and the bottom of the container, whereas the top of the container is assumed to be open.
The fluid velocity at the ring surface corresponds to the  structural velocity.
In terms of the mesh deformation problem the mesh nodes on the container and walls of the phantom domains are restricted to a vertical movement.
The structural deformation is prescribed as a Dirichlet value for the ring boundary.\\
\\
The FSI problem is solved in a partitioned solution approach \cite{FelippaParkFarhat1998}.
On the structural side, the deformation of the ring are represented by a linear elastic problem solved with isogeometric analysis (IGA) \cite{HughesCottrellBazilevs2004}.
The contact interaction between the ring and the bottom of the container is considered via the penalty method \cite{TemizerWriggersHughes2011}.
The flow field induced by the motion of the ring is described by the Navier-Stokes equations which are solved by the DSD/SST approach in combination with the presented PD-DMUM.
The two field problems are strongly coupled in time \cite{Wall1999}, and for the spatial coupling we apply a NURBS-based coupling following \cite{HostersEtAl2017}.\\
\\
In Figures \ref{fig:Ball1} to \ref{fig:Ball5}, we present snapshots of the simulation at different points in time, starting from the initial position of the ring, via the moment when the ring is in contact with the bottom of the container, up to the point of maximal altitude after the first contact interaction.
As it can be guessed from the snapshot in Figure \ref{fig:Ball3}, one element remains between the bottom of the container and the falling ring.
This element will not be removed because we cannot exactly comply with the contact conditions using the penalty method.
Nevertheless, it can be observed in every snapshot, that mesh cells experience large displacements but only little deformations.
Due to the application of the PD-DMUM, the entire FSI problem was solved without remeshing.
\begin{figure}[h!]
	\centering
	\begin{minipage}{.4\textwidth}
		\centering
		\includegraphics[width=\textwidth]{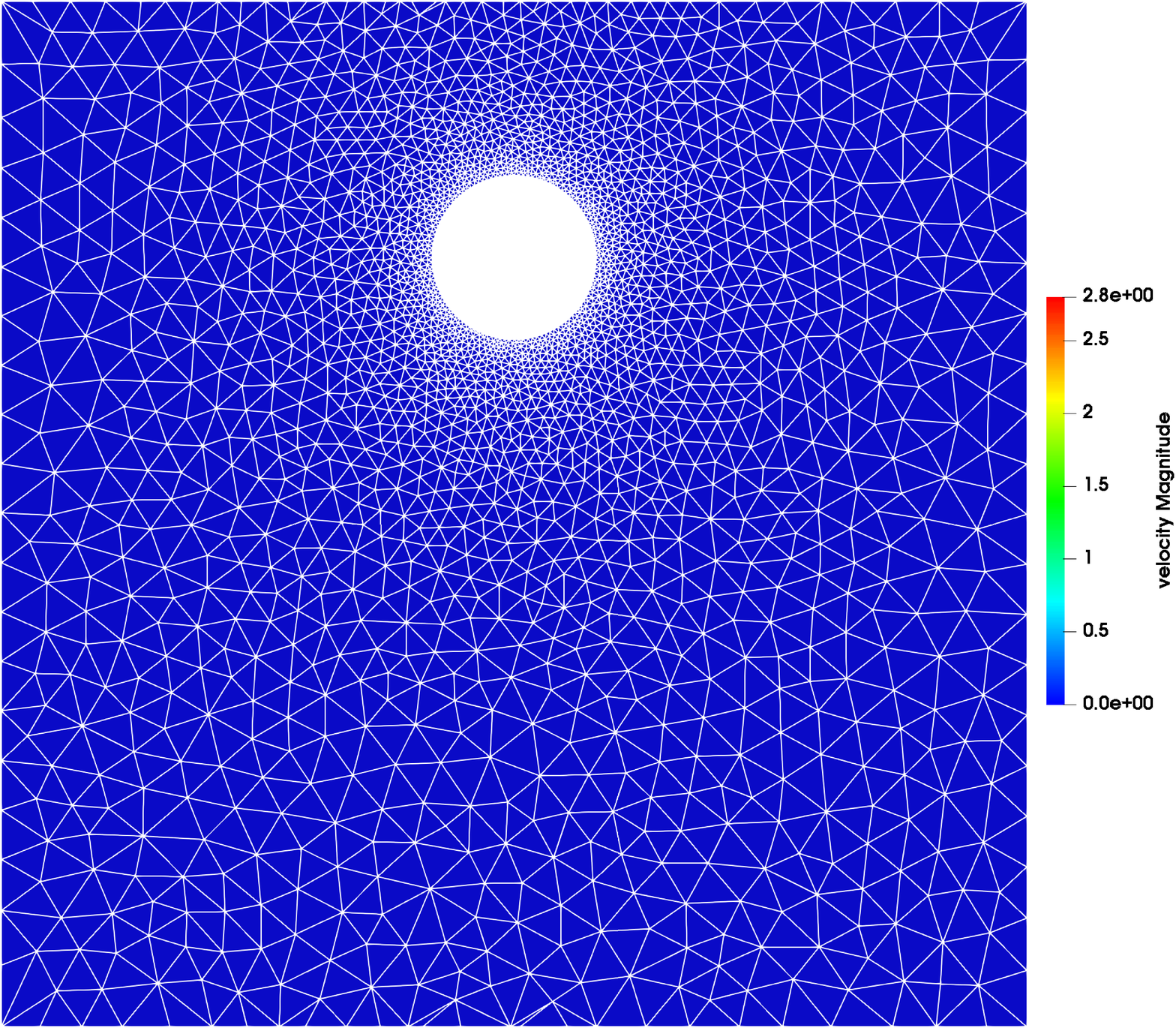}
		\caption{Velocity at $t=0 ~\text{s}$.}
		\label{fig:Ball1}
	\end{minipage}%
	\begin{minipage}{.4\textwidth}
		\centering
		\includegraphics[width=\textwidth]{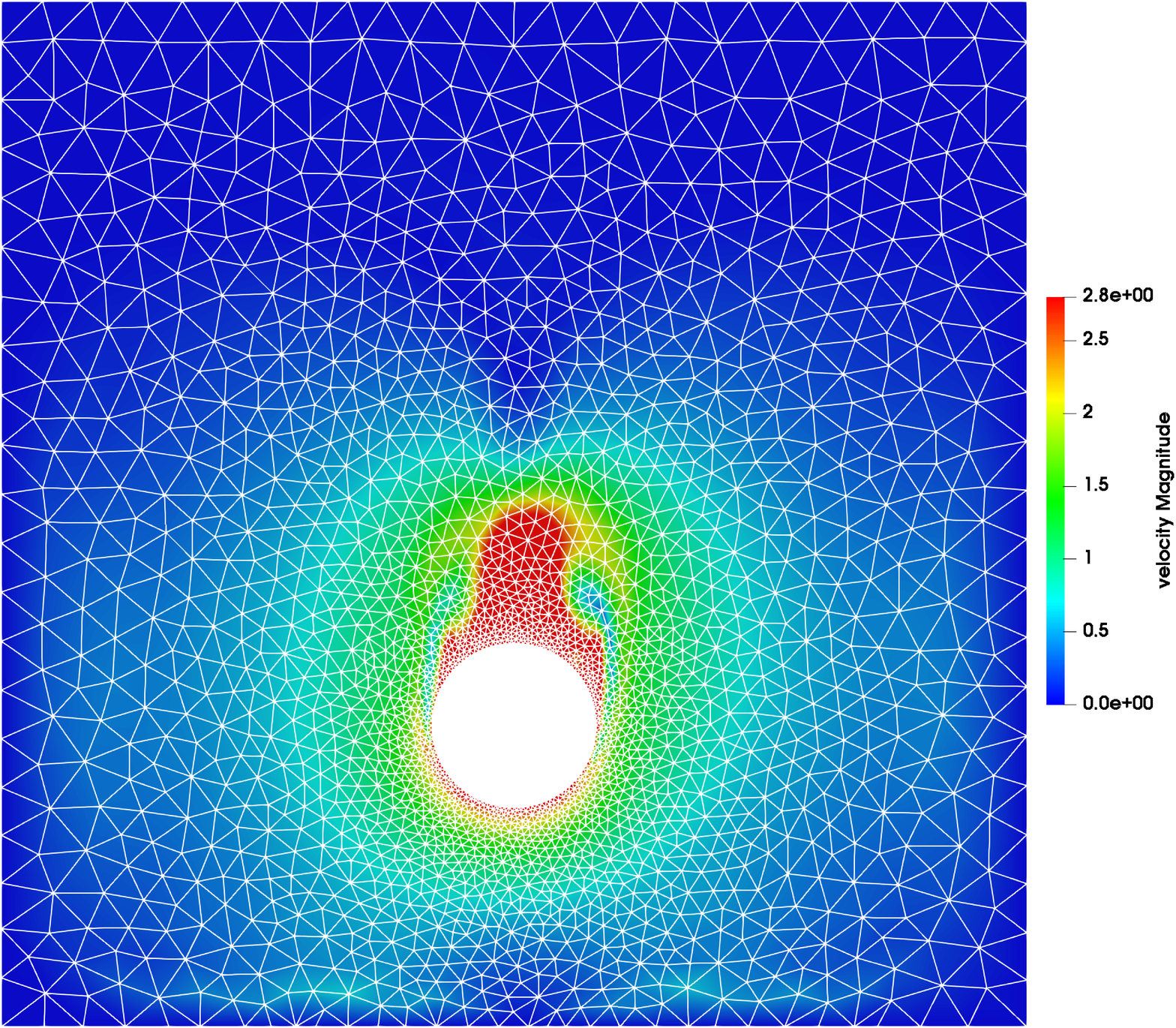}
		\caption{Velocity at $t=0.55~\text{s}$.}
		\label{fig:Ball2}
	\end{minipage}
	\\
	\begin{minipage}{.4\textwidth}
		\centering
		\includegraphics[width=\textwidth]{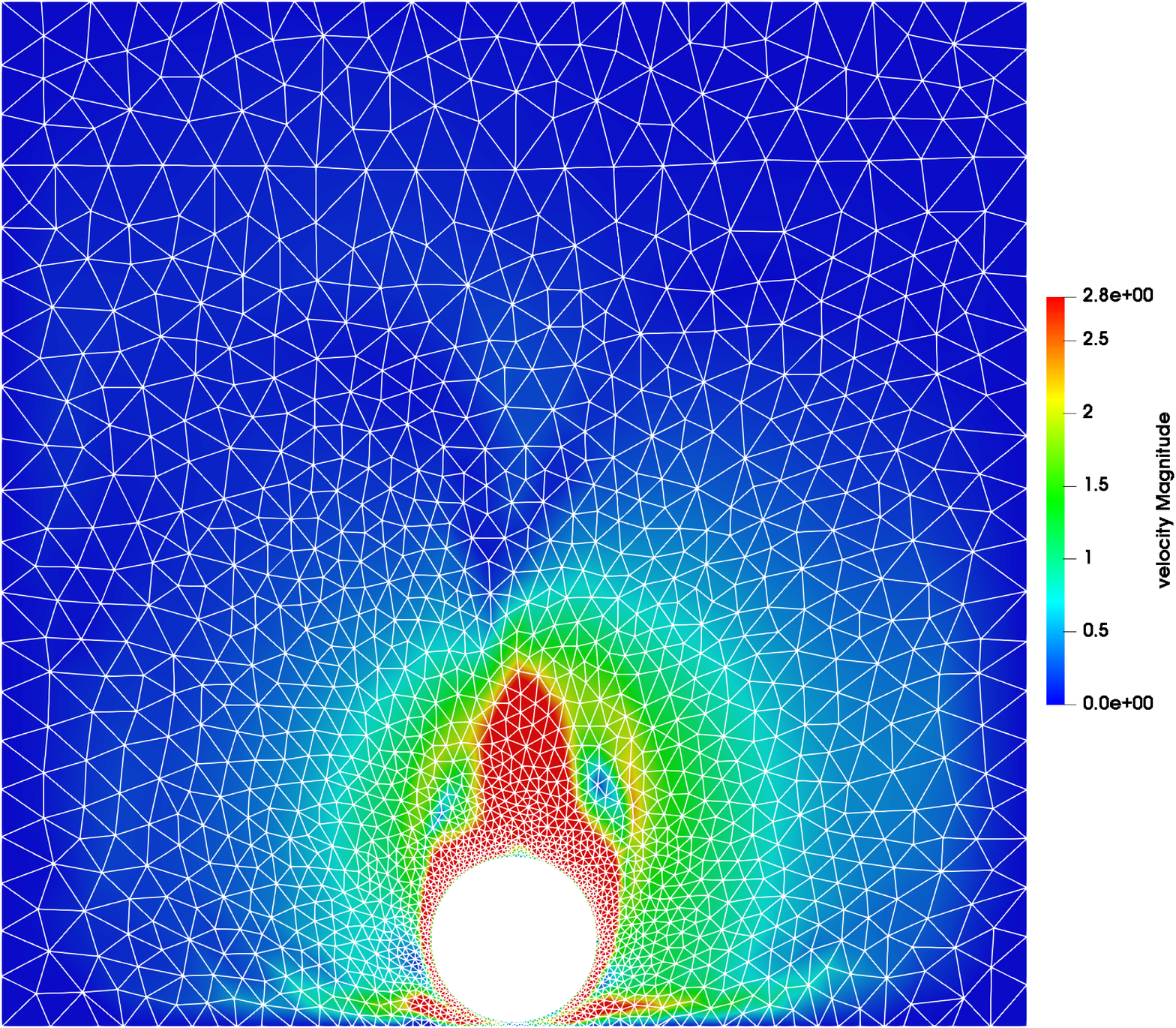}
		\caption{Velocity at $t=0.75~ \text{s}$.}
		\label{fig:Ball3}
	\end{minipage}%
	\begin{minipage}{.4\textwidth}
		\centering
		\includegraphics[width=\textwidth]{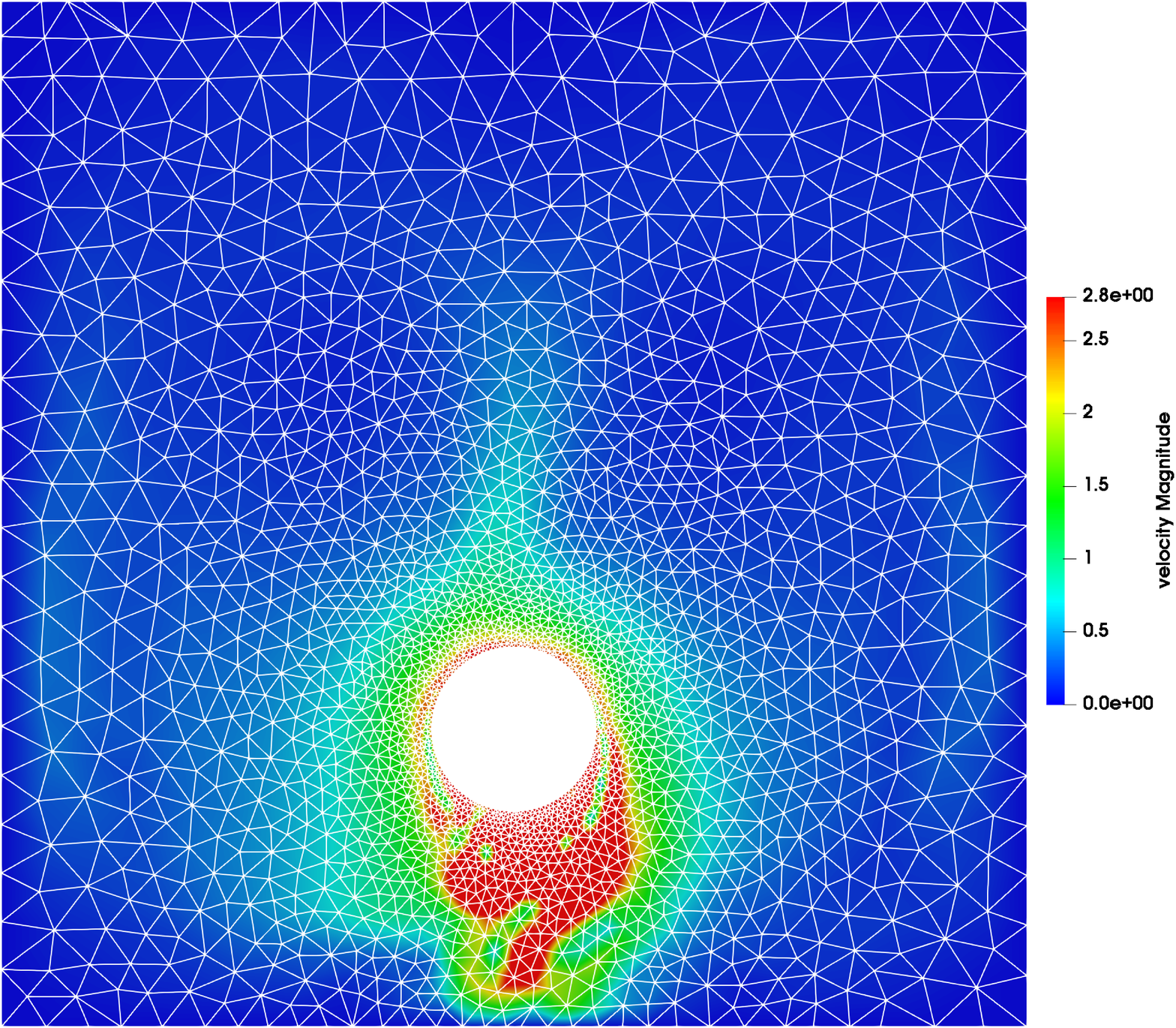}
		\caption{Velocity at $t=1.0 ~\text{s}$.}
		\label{fig:Ball4}
	\end{minipage}\\
	\centering
\end{figure}
\begin{figure}[h!]
	\centering
	\ContinuedFloat
	\begin{minipage}{.4\textwidth}
		\centering
		\includegraphics[width=\textwidth]{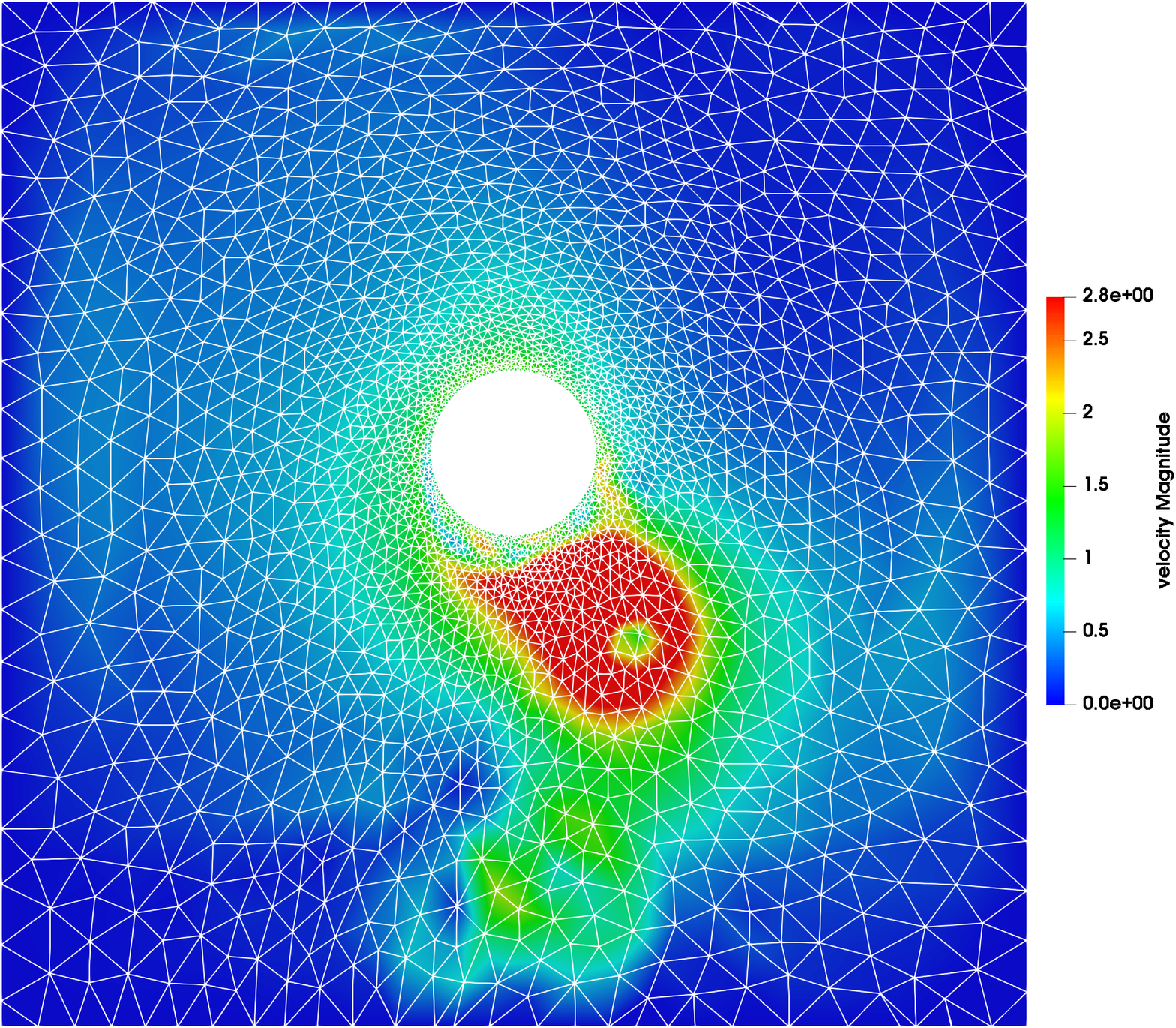}
		\caption{Velocity at $t=1.45~\text{s}$.}
		\label{fig:Ball5}
	\end{minipage}
\end{figure}
\\

\section{Discussion}
\label{sec:5}
In this paper, we presented a novel approach for the mesh update of boundary conforming meshes, particularly developed for  problems with large unidirectional boundary movements, the PD-DMUM.
Subsequent to the description, we evaluated the PD-DMUM in two test cases.
In the first test case we showed by means of a Poiseuille flow the general agreement of the PD-DMUM with results of consisting methods.
In the fluid-structure interaction problem presented in the second test case we emphasised the applicability of the PD-DMUM in complex processes with moving boundaries.\\
\\
Both test cases yielded consistently good results.
Although, the PD-DMUM still requires interpolation at single element nodes, it meets its two major challenges: (1) Even in complex processes as structural contact, remeshing of the domain is entirely avoided, and (2) the computed solution of the flow problems is in accordance with solutions computed by conventional mesh update methods.\\
\\
In summary we successfully introduced a new mesh update approach, where the first test cases congruently showed good results.
Further this method bears potential for problem specific improvements of the computational effort, by combination of the  PD-DMUM with the virtual ring or other mesh update approaches that reduce the deforming mesh area.

\begin{acknowledgement}
This work was supported by the German Research Foundation under the Cluster of Excellence "Integrative production technology for high-wage countries" (EXC128) as well as the  German Research Foundation under the Cluster of Excellence "Internet of Production". Computing resources were provided by the AICES graduate school and RWTH Aachen University Center for Computing.
\end{acknowledgement}

\end{document}